\documentclass [11pt] {article}

\usepackage{amsfonts}
\usepackage{amsmath}
\usepackage{amssymb}
\usepackage{color}
\usepackage{graphicx}
\usepackage{times}
\usepackage[round]{natbib}
\usepackage[T1]{fontenc}
\usepackage{url}
\usepackage{endnotes}

\begin{document}

\title{SIAM 2020\footnote {Submitted to \textit {SIAM News}.}}

\author{Joseph F. Grcar\,\footnote{6059 Castlebrook Drive, Castro Valley, CA 94552 USA.}\, \footnote {\texttt {jfgrcar@comcast.net}, \texttt {jfgrcar@gmail.com.}}}

\date {}

\maketitle

\begin{abstract}
\noindent
Some observations are made about how the Society for Industrial and Applied Mathematics (SIAM) might be better oriented to serve the industrial and interdisciplinary mathematics community in the future.


\end{abstract}


\paragraph {Introduction} When Bill Coughran and Gene Golub planned the SIAM 30th anniversary meeting at Stanford, I suggested to Gene that SIAM conferences would be more interesting if members could organize their own sessions. Gene took my advice and the name that I suggested for the sessions, ``minisymposia.'' The Stanford meeting in 1982 was the first to have them. The only other change at SIAM that benefited members so directly was the creation of the activity groups which occurred around the same time. 

Now is a good time to increase the pace of innovation. Academic members know that in recent decades the number of upper division mathematics students has decreased so much that one research university proposed closing its mathematics department. David Bressoud,\endnote {D. M. Bressoud, ``Is the Sky Still Falling?,'' \textit {Notices of the AMS}, 2009, 56:1 20--25.} the president of the mathematics association, recently noted that ``the situation is far from healthy, and in many respects we are worse off now than we were in 1995.'' SIAM members on the other hand teach the interesting courses and most of them are not even in mathematics departments. 

The future growth of industrial mathematics lies in interdisciplinary studies which SIAM members are best suited to address by inclination and situation.  The sociology of science shows that academic fields result from a demand for teaching specialists and from a professional ethic that researchers share across institutional boundaries.  SIAM needs to evolve to keep up with its members as the interdisciplinary subjects mature.  For example, in today's economy the ``I'' in SIAM might just as well stand for the information industries whose success depends on mathematics broadly interpreted.  With the winds of change blowing in SIAM's direction I am going to propose a vision for SIAM to surge ahead in the 21st century.

My vision for SIAM has four points: education, industry, leadership, and publishing. I think it will take a decade for the changes to attain their potential, so I am going to call my vision SIAM 2020. I propose to build on the features of the present SIAM and make them more effective.

\newcommand {\edu} [1] {(\textbf {E.#1})}
\newcommand {\ind} [1] {(\textbf {I.#1})}
\newcommand {\lea} [1] {(\textbf {L.#1})}
\newcommand {\pub} [1] {(\textbf {P.#1})}

\paragraph {Education} Lynn Steen is a former mathematics association president who has written with some astonishment of the exciting mathematics ``stealth curriculum that thrives outside the confining boundaries of college and university mathematics departments.''\endnote {L.\ A.\ Steen, ``Facing Facts: Achieving Balance in High School Mathematics,'' \textit {Mathematics Teacher\/}, Dec-Jan 2005-2006, 100th anniversary issue, 86--95. Steen's quotation paraphrases another article by him on university education.} Lynn refers to the kind of mathematics known to SIAM members but ``stealthy'' to many because it has no focal point.  Much of industrial mathematics is inherently interdisciplinary so it currently has no advocates in agencies like the National Science Foundation and is only beginning to be embodied in formal study programs outside university mathematics departments. 

SIAM should address the policy aspects of mathematics education and represent its community in curriculum matters. I propose a three-step plan. \edu 1 An activity group should be formed for industrial and interdisciplinary mathematics education (IIME). Some of the  most heavily attended sessions at the New Orleans annual meeting were about computational science and engineering education, so members evidence strong interest in education when given a venue. Many SIAM members can gather the necessary signatures just by walking down the hall. 

\edu 2 SIAM should convey the views of industrial and interdisciplinary mathematicians at any national meetings that discuss mathematics curricula. For example, participants at a recent meeting sponsored by the Association of American Universities were astonished to find that natural and social scientists disagreed with leaders from traditional mathematics departments over the curriculum for high school math.\endnote {D.\ T.\ Conley, ed., \textit {Understanding University Success, A project of the Association of American Universities and The Pew Charitable Trusts,\/} Center for Educational Policy Research, University of Oregon, 2003.} 

\edu 3 SIAM should become the accreditation service for interdisciplinary mathematics programs. Many universities stress these programs but currently they have no accreditation even though they are usually associated with professional colleges of engineering. These programs will be accredited by someone eventually, so it is in everyone's best interests that SIAM should begin accrediting them now. This service could be done in cooperation with consortia such as ABET and CSAB that accredit many university departments including those at the universities where most SIAM members are employed.

\paragraph {Industry} Some members may be surprised to learn that the industrial aspect of SIAM has two components: members who work in industry (12\% of non-student members), and academic members in non-mathematics departments (28\%) which typically have strong ties to industry.\endnote {Completing the membership data, of approximately 7900 non-student members, 38\% are in mathematics departments, 11\% are outside academe and industry such as in government or government laboratories, and 11\% do not declare employers.} I suspect there is a third class of potential members who studied mathematics but do not use it directly in their professions. Growing all these classes of members depends on making SIAM more prominent in STEM\endnote {STEM is higher education jargon for science, technology, engineering, and mathematics.} education and on retaining members who may not contribute to formal mathematics journals.

I propose two ways to make SIAM more attractive to such members. \ind 1 \textit {SIAM News\/} should become a magazine of general interest for industrial and interdisciplinary mathematics such as computational science. The excellent material that appears there now would have greater impact and would attract a larger volume of like material in a citable magazine format that is easily archived and displayed on the internet. Many societies have such magazines for articles, commentary, news, reports, and reviews. 

\ind 2 The SIAM fellows program should honor members who are industry leaders. There is room to grow because 12\% of non-student members work in industry while only 6\% of fellows have industrial ties. The fellowship committee only acts on nominations it receives, so it should fall to the officers and to the SIAM vice president for industry to see that industrial leaders are appropriately honored by arranging winning nominating packages for them. A magazine of interest to those who know industrial mathematics but who might not work with formal mathematics, and a fellows program recognizing their lifetime contributions, would make mathematics more prominent generally.

\paragraph {Leadership} University mathematics departments with few exceptions shun mathematics that is closely related to other fields.\endnote {The following institution is chosen for example only because its web pages are so well constructed to permit searches.  The University of Texas at Austin taught 7 courses on finite element analysis in the fall 2010 semester. None were offered in the mathematics department.} Whether annoying or amusing, this process helps demarcate academic disciplines.\endnote {C.\ Calhoun, ``The specificity of American higher education,'' in R.\ Kalleberg, et al., \textit {Comparative Perspectives on Universities\/}, JAI Press, Stamford, 2000, 47--82.} Consequently many branches of mathematics are allied with other fields and are represented through interest groups in larger societies. These subjects don't apply mathematics, they create it, and their contribution of applicable mathematics should be recognized when setting national policy for education and research.

The other mathematics societies are unlikely to pursue this course, so I propose SIAM should assert leadership, in three steps. \lea 1 SIAM should explore organizing a conference board of industrial and interdisciplinary mathematics societies and interest groups to identify common concerns and to attest widespread research activities in mathematics outside mathematics departments. Indeed, many of the topics of research are embodied in the titles of SIAM's own journals and activity groups. 

\lea 2 SIAM should lead in obtaining the appointment of an industrial mathematician to the National Science Board, which is the governing body for the National Science Foundation. SIAM is the appropriate society to make the nomination because of its interdisciplinary credentials; the other mathematics societies should be ready to endorse a nomination because no mathematician serves on the board. 

\lea 3 SIAM should be more transparent in its own governance. Since it has members from many research communities, SIAM officers tend to look for participation from insiders like themselves and from the groups they know, which has led many people to characterize SIAM as cliquish. The SIAM president recently issued a broad call for nominations to elected positions. This call needs to be extended to the many more appointed positions both on committees and editorial boards whose selection process is not transparent. 


\paragraph {Publishing} Many professional societies devote considerable effort to two dissimilar businesses: conference management and technical publishing. Their publishing business is threatened two ways. The low cost of internet distribution facilitates entry of competing outlets for papers, and the decreasing budgets of universities limits the customers for journals. Societies charge relatively high fees for access to comparatively small collections.\endnote {ScienceDirect charges 4-year colleges \$36,470 annually for internet access to 1860 journals, while SIAM charges \$5,928 for access to 15 journals.} For example, SIAM has only one journal in the top 5 by volume of articles about any single branch of mathematics.\endnote {\textit {SINUM\/} publishes 2.3\% of all papers whose primary mathematics subject classification is 65, numerical analysis, compared to \textit {Applied Mathematics and Computation\/} which publishes 9.9\%.} The only way to survive is by increasing the quantity of published papers, thereby becoming an obvious choice for submission and an indispensable choice for continued purchase. 

I propose the twin goals that within the next decade SIAM should be in the top 5 publishers of research papers about each branch of industrial mathematics and it should be the dominant publisher across branches.  I propose two steps to reach these goals.  The first step will be described in some length.

\pub 1 Industrial and interdisciplinary mathematics is difficult to address because it is fragmented into many specialities. SIAM publishing, however, has unique advantages in the meetings and membership of the professional side of the society. Authors want their papers to quickly contribute to their C.V.'s where each paper can potentially add two items: one for printing the paper and one for presenting it. SIAM currently does not offer both because the publishing business is separate from the conference business. I propose linking the two in a special class of conference papers that are peer-reviewed on an expedited schedule for both publication and presentation. Only SIAM can do this because it has the breadth of expertise among its members for the reviewing and the breadth of conferences for the presenting. Each activity group should gather reviewers who promise rapid evaluations after which the accepted papers appear in \textit {SIAM Lecture Notes\/} on the topic of the activity group. For example, \textit {SIAM Lecture Notes on Computational Science and Engineering\/} could become hallmarks of the society. The lecture notes can be sold either individually or in subscription. 

This business model is the one that Springer Verlag uses to dominate computer science publishing through the ubiquitous \textit {Lecture Notes in Computer Science}. Some of them are monographs but many are proceedings for recurring conferences that are recognized as high-quality publications (Springer does not organize the conferences). Only a minority of SIAM members are employed in mathematics departments where traditional journals are the most highly regarded venues; most SIAM members belong to engineering departments or to industry where peer-reviewed conference publications are equally regarded.

This plan has the added benefit of offering professional development to more members. The officers of activity groups can learn editorial and publishing skills by managing the Lecture Notes. In this way the pool of experienced editors for the journals is also enlarged.

\pub 2 SIAM journals have not done well in the long run because many were the first in their subjects and yet they did not grow with their fields. Anecdotal evidence suggests the review process is too long and the editorial tradition does not recognize the interests of authors.\endnote {For discussion of best editorial practice see (a) M.\ L.\ Cooper, ``Problems, Pitfalls, and Promise in the Peer-Review Process,'' \textit {Perspectives on Psychological Science\/}, 2009, 4:1 84-90, and (b) E.\ W.\ K.\ Tsang and B. S. Frey, ``The As-Is Journal Review Process: Let Authors Own Their Ideas,'' \textit {Academy of Management Learning \& Education\/}, 2007, 6:1 128--136.} I suggest that all submissions should be double blind (so reviewers and associate editors do not know the authors) and the criteria for acceptance or revision should be made explicit in an author's bill of rights.

\paragraph {Conclusion} In summary, strong academic disciplines result from a demand for specialists to teach the university curriculum and an extramural professional society that represents the collective views of the discipline across universities. With SIAM 2020 both can thrive in the next century.

\theendnotes

\end{document}